# GAUSSIAN PROCESSES, KINEMATIC FORMULAE AND POINCARÉ'S LIMIT

BY JONATHAN E. TAYLOR[1,2] AND ROBERT J. ADLER[1]

*Stanford University and Technion*

We consider vector valued, unit variance Gaussian processes defined over stratified manifolds and the geometry of their excursion sets. In particular, we develop an explicit formula for the expectation of all the Lipschitz–Killing curvatures of these sets. Whereas our motivation is primarily probabilistic, with statistical applications in the background, this formula has also an interpretation as a version of the classic kinematic fundamental formula of integral geometry. All of these aspects are developed in the paper.

Particularly novel is the method of proof, which is based on a an approximation to the canonical Gaussian process on the $n$-sphere. The $n \to \infty$ limit, which gives the final result, is handled via recent extensions of the classic Poincaré limit theorem.

**1. Introduction.** The central aim of this paper is to describe a new result in the theory of Gaussian related fields, along with some of its implications to both geometry, and to a lesser extent, to statistics.

The basic object of interest is a $\mathbb{R}^k$ valued random field $y$ defined on a parameter space $M$ and its *excursion sets*

$$(1.1) \qquad A(f, M, D) \stackrel{\Delta}{=} \{t \in M : y(t) \in D\},$$

where $D \subset \mathbb{R}^k$. For most of the paper, we shall take $M$ and $D$ to be $C^2$ stratified manifolds in $\mathbb{R}^N$ and $\mathbb{R}^k$, respectively.

---

Received January 2008; revised May 2008.
[1]Supported in part by US–Israel Binational Science Foundation Grant 2004064.
[2]Supported in part by NSF Grant DMS-04-05970, and the Natural Sciences and Engineering Research Council of Canada.
*AMS 2000 subject classifications.* Primary 60G15, 60G60, 53A17, 58A05; secondary 60G17, 62M40, 60G70.
*Key words and phrases.* Gaussian fields, kinematic formulae, excursion sets, Poincaré's limit, Euler characteristic, intrinsic volumes, geometry.

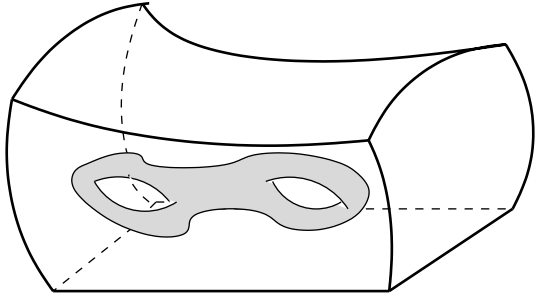

Fig. 1.  *A saggy couch under stress: A stratified manifold with its excursion sets.*

Stratified manifolds are basically sets that can be partitioned into the disjoint union of manifolds, so that we can write

$$(1.2) \qquad M = \bigsqcup_{j=0}^{\dim M} \partial_j M,$$

where each stratum, $\partial_j M$, $0 \le j \le \dim(M)$, is itself a disjoint union of a number of $j$-dimensional manifolds. A typical 3-dimensional example is given by the saggy couch of Figure 1, in which case $\partial_3 M$ is the interior of the couch; $\partial_2 M$ the collection of the six sides, some concave and some convex; $\partial_1 M$ is made up of the 12 edges; and $\partial_0 M$ contains the 8 corner vertices. In Figure 1, an excursion set might be the grey area where the stress $y$ is greatest.

Our aim is to study the global geometry of excursion sets, as measured through their Lipschitz–Killing curvatures, $\mathcal{L}_j(A(f,M;D))$, $j = 0,\ldots,\dim(M)$. In particular, since these curvatures are random variables, we shall be interested in computing their expectations. We shall define Lipschitz–Killing curvatures below. However, if you are unfamiliar with them, at this stage it suffices to know that $\mathcal{L}_N(A)$ is a measure of the volume of $A$, $\mathcal{L}_{N-1}(A)$ a measure of its surface area, and $\mathcal{L}_0(A)$ its Euler characteristic, an important topological invariant.

We cannot do this for all $y$. For a start, $y$ will both have to be smooth enough for basic differential geometric techniques to be applicable. Thus, a basic requirement will be that $y$ has, with probability one, $C^2$ sample paths. Furthermore, writing $y = (y_1,\ldots,y_k)$, we shall assume that the $y_i$ are independent, identically distributed (hereafter i.i.d.) centered Gaussian processes of constant variance, which we take to be 1. The processes are, however, not assumed to be stationary. For such a $y$, we shall prove that

$$(1.3) \quad \mathbb{E}\{\mathcal{L}_i(A(y,M,D))\} = \sum_{j=0}^{\dim M - i} \begin{bmatrix} i+j \\ j \end{bmatrix} (2\pi)^{-j/2} \mathcal{L}_{i+j}(M) \mathcal{M}_j^\gamma(D),$$

where the combinatorial flag coefficients are defined below at (4.3) and the $\mathcal{M}_j^\gamma$, described and defined in Sections 3 and 6.1, are certain (Gaussian)

Minkowski functionals that to a certain extent, play the role of Lipschitz–Killing curvatures in Gauss space. The Lipschitz–Killing curvatures on both sides of (1.3) are computed with respect to a specific Riemannian metric induced on $M$ by the component processes $y_j$. Note, however, that $\mathcal{L}_0(A)$ is the Euler–Poincaré characteristic of $A$, and so independent of any Riemannian structure [cf. Theorem 4.1 for a formal statement of (1.3)].

1.1. *What is new here?* The result (1.3) has a long history. If $M$ is a simple interval $[0, T]$, $y$ is real valued and stationary, and $D = [u, \infty)$, then (1.3) is essentially the famous Rice formula, which gives the mean number of upcrossings of the level $u$ by $f$, and dates back to 1939 [15] and 1945 [16]. Since then, there have been tens, if not hundreds, of papers extending the original Rice formula in many ways, with the developments up until 1980 summarized in [1]. More recently, there was a series of papers by Worsley (e.g., [24, 25, 26, 28]) that were important precursors to the general theory presented in this paper. However, the first precursor to (1.3), at the level of processes over manifolds with $C^2$ boundaries, appeared only in 2002 in [21], where we considered only the first Lipschitz–Killing curvature $\mathcal{L}_0(A_u(f, M))$ and then only for real valued $y$. In [20], one of us (JET) extended this to vector valued $y$, which allowed for the derivation of the far more general, and far more elegant, (1.3) for the first Lipschitz–Killing curvature.

What is new here then is the extension to parameter spaces as general as stratified manifolds, and the extension to all Lipschitz–Killing curvatures. Both of these are important for applications. However, perhaps more important, and certainly more novel than either of these, is the method of proof. The proofs in the current paper are new, and far more geometric than the earlier ones. In particular, the proof in [20] progressed primarily by evaluating both sides of (1.3) and then showing that they were equivalent. The current proof starts on the left-hand side and, eventually, yields the right-hand side. The geometric nature of the current proof also explains *why* the two sides *should* be equal.

1.2. *Statistical implications.* The general structure of (1.3) has significant implications for a class of problems out of the purely Gaussian scenario. Taking $F : \mathbb{R}^k \to \mathbb{R}$ to be piecewise $C^2$, along with appropriate side conditions, and defining a (now non-Gaussian) process

$$(1.4) \qquad f(t) = F(y(t)) = F(y_1(t), \ldots, y_k(t)),$$

with $y$ Gaussian as above, it follows immediately from (1.3) that

$$(1.5) \quad \begin{aligned} &\mathbb{E}\{\mathcal{L}_i(A(f, M, [u, \infty)))\} \\ &\quad = \sum_{j=0}^{\dim M - i} \begin{bmatrix} i+j \\ j \end{bmatrix} (2\pi)^{-j/2} \mathcal{L}_{i+j}(M) \mathcal{M}_j^\gamma(F^{-1}[u, +\infty)). \end{aligned}$$

Non-Gaussian processes of the form (1.4) appear naturally in a wide variety of statistical applications of smooth random fields (e.g., [2, 3, 4, 19, 24, 25, 26] with an excellent introductory review in [27]).

An additional and extremely important application of (1.3) lies in the so called "Euler characteristic heuristic" that for a wide range of random fields $f$,

$$\left|\mathbb{P}\Big\{\sup_{t\in M} f(t) \geq u\Big\} - \mathbb{E}\{\mathcal{L}_0(A(f, M, [u, \infty)))\}\right| \leq error(u),$$

where $error(u)$ of a smaller order than both of the other terms as $u \to \infty$. In the Gaussian case, this heuristic is now a well-established theorem, and the error term is known to be of order $\exp(-u^2(1+\eta)/2)$ (for an identifiable $\eta > 0$) while both the probability and expectation are of order $\exp(-u^2/2)$ [22]. The ability to compute the expectation therefore provides useful, explicit approximations for the excursion probability.

1.3. *Geometry.* One of the basic results of integral geometry is the so-called kinematic fundamental formula (henceforth KFF), which in its simplest form, states that for nice subsets $M_1$ and $M_2$ of $\mathbb{R}^n$,

$$
\begin{aligned}
(1.6) \quad &\int_{G_n} \mathcal{L}_i(M_1 \cap g_n M_2) \, d\nu_n(g_n) \\
&= \sum_{j=0}^{n-i} \begin{bmatrix} i+j \\ i \end{bmatrix} \begin{bmatrix} n \\ j \end{bmatrix}^{-1} \mathcal{L}_{i+j}(M_1)\mathcal{L}_{n-j}(M_2).
\end{aligned}
$$

Here, $G_n$ is the isometry group of $\mathbb{R}^n$ with Haar measure $\nu_n$ normalized so that for any $x \in \mathbb{R}^n$ and any Borel $A \subset \mathbb{R}^n$, $\nu_n(\{g_n \in G_n : g_n x \in A\}) = \mathcal{H}_n(A)$, where $\mathcal{H}_n$ is $n$-dimensional Hausdorff measure. (See [12, 18] for $M_j$ elements of the convex ring or similar, and [6] for more esoteric $M_j$ closer to the spirit of this paper.)

Now reconsider (1.3). Taking $(\Omega, \mathcal{F}, \mathbb{P})$ as the probability space on which $y$ lives, (1.3) can be rewritten as

$$
\begin{aligned}
(1.7) \quad &\int_\Omega \mathcal{L}_i(M \cap (y(\omega))^{-1} D) \, d\mathbb{P}(\omega) \\
&= \sum_{j=0}^{\dim M - i} \begin{bmatrix} i+j \\ j \end{bmatrix} (2\pi)^{-j/2} \mathcal{L}_{i+j}(M) \mathcal{M}_j^\gamma(D).
\end{aligned}
$$

Written this way, it is clear on comparing (1.6) and (1.7) that our main result can now be interpreted as a KFF over Gaussian function space, rather than over the isometry group on Euclidean space. We find this interpretation novel and intriguing bridging as it does between a probabilistic problem and a geometric answer of classic form.

1.4. *More on stratified manifolds.* Although (1.2) gives the basic structure of a stratified manifold, certain assumptions need to be made for the results of this paper to hold. In particular, we need to assume that each connected component in each stratum is a $C^2$ manifold. More importantly, we need to assume that both $M$ and $D$ are Whitney stratified spaces (see [10, 14]) which ensures that the various strata are "glued together" in a smooth way. There are further technical assumptions of "tameness" and "regularity" that can be found in Chapter 15 of [3]. One of these is local convexity; a demand much weaker than convexity (the saggy couch is locally convex), which while not necessary, is something that we shall assume in this paper so as to make many formulae somewhat more manageable. Handling all these technicalities is not only tiresome, consuming both space and time, but often also quite difficult. Nevertheless, we shall completely avoid them since our main aim in the current paper is primarily didactic. We want to present the main result in a setting that gives the flavor of the general result, introduces the new style of proof and yet does so in a way that is both self-contained and comparatively easy to follow. If you want all the details, then Chapter 15 of [3] is the place to turn. If you want rigor without leaving this paper, then you can assume that $M$ is a manifold without boundary (necessarily locally convex) such as a sphere or a simple cube, which has all the edges and corners but none of the curvature of Figure 1. In neither of these cases are additional assumptions needed, and for both of them the results and proofs of the paper are still new.

1.5. *A plan of action.* Our first task, in the following two sections, will be to define the Lipschitz–Killing curvatures $\mathcal{L}_j$ and the Gaussian Minkowski functionals $\mathcal{M}_j^\gamma$. Once this is done, then at very least all the terms in (1.3) will be clearly defined, and we shall be able to state the formal result as Theorem 4.1. Thus, the reader who is not at all interested in proofs can leave us at that point.

We shall then slowly start on the proof, proving in detail only a special case of the Theorem 4.1 for the so-called isotropic process on the unit sphere in $\mathbb{R}^N$. This is the most interesting part of the paper, and contains not only new results, but more importantly, a conceptually new way of looking at things. Despite limiting the proof to one special case, we shall explain in Section 8 why this case is particularly important, and also why all other cases can, in principle (but with quite a lot of additional work) be shown to follow from it.

However, even for this specific choice of process, the computations are far from trivial and so we approach them in an indirect fashion. In particular, we shall start by looking at non-Gaussian processes with finite expansions and with coefficients coming from random variables distributed uniformly over high dimensional spheres. This will enable us to take expectations using an appropriate version of the KFF on spheres.

The passage from this scenario to the Gaussian one will be via a limit theorem for projections of uniform random variables on the $n$-sphere as $n \to \infty$, historically associated with the name of Poincaré, although we shall need a slight extension of some more recent versions due to Diaconis and Freedman [7] and their generalizations to matrices in [8].

**2. Lipschitz–Killing curvatures.** Perhaps the easiest way to meet Lipschitz–Killing curvatures is via Weyl's tube formula, which in our setting states that for quite wide classes of sets $M$ in $\mathbb{R}^l$, and for sufficiently small $\rho \geq 0$, there exist numbers $\mathcal{L}_j(M)$ for which

$$(2.1) \qquad \mathcal{H}_l(\mathrm{Tube}(M,\rho)) = \sum_{i=0}^{\dim M} \rho^{l-i} \omega_{l-i} \mathcal{L}_i(M),$$

where $\mathrm{Tube}(M,\rho) \triangleq \{t \in \mathbb{R}^l : \inf_{s \in M} |s-t| \leq \rho\}$, and $\omega_m$ is the volume of the unit ball in $\mathbb{R}^m$ [cf. (4.3)] and, in a notation anticipating something more general, $\mathcal{H}_l$ is $l$-dimensional Lebesgue measure.

Weyl's formula actually defines the $\mathcal{L}_j(M)$ for this Euclidean case for all $0 \leq j \leq \dim(M)$, and we define $\mathcal{L}_j(M) \equiv 0$ for all $j > \dim(M)$. However, they can also be computed directly, and have a natural extension to the situation in which $M$ is a Riemannian manifold, a situation which we shall soon need, even though our parameter spaces $M$ are basically Euclidean. In the Riemannian case, we shall need to assume that $M$ is embedded in an ambient manifold $\widetilde{M}$. We shall denote the Riemannian metrics by $g$ and $\widetilde{g}$, respectively, the curvature tensors of $M$ and $\widetilde{M}$ by $R$ and $\widetilde{R}$, and denote the scalar second fundamental forms of the strata $\partial_j M$ as they sit in $\widetilde{M}$ by $\widetilde{S}$. We also have Riemannian (Hausdorff) volumes on $M$ and $\partial_j M$, which we continue to denote by $\mathcal{H}_N$ and $\mathcal{H}_j$, respectively.

In this more general setting, Lipschitz–Killing curvatures of $M$ are defined by

$$\begin{aligned}
\mathcal{L}_i(M) = \sum_{j=i}^{\dim M} (2\pi)^{-(j-i)/2} \\
\times \sum_{m=0}^{\lfloor (j-i)/2 \rfloor} \frac{(-1)^m C(l-j, j-i-2m)}{m!(j-i-2m)!} \\
\times \int_{\partial_j M} \int_{S(T_t \partial_j M^\perp)} \mathrm{Tr}^{T_t \partial_j M}(\widetilde{R}^m \widetilde{S}_\eta^{j-i-2m}) \\
\times \mathbb{1}_{N_t M} \mathcal{H}_{l-j-1}(d\eta) \mathcal{H}_j(dt),
\end{aligned} \qquad (2.2)$$

where $T_t\,\partial_j M$ is the tangent space to $\partial_j M$ at $t$, $T_t\,\partial_j M^\perp$ is the normal space, $N_t M$ is the normal cone to $M$ at $t$, $\text{Tr}^{T_t\,\partial_j M}$ is the trace in $T_t\,\partial_j M$, and $S(T_t\,\partial_j M^\perp)$ is the unit sphere in $T_t\,\partial_j M^\perp$. [We shall generally write $S_\lambda(V)$ and $B_\lambda(V)$ for the sphere and ball of radius $\lambda$ in a vector space $V$, with $S \equiv S_1$ and $B \equiv B_1$.] The constants $C(m,i)$ are given by

$$C(m,i) \triangleq \begin{cases} \dfrac{(2\pi)^{i/2}}{s_{m+i}}, & m+i > 0, \\ 1, & m = 0, \end{cases}$$

and $s_n = \Gamma(\frac{n}{2})/2\pi^{n/2}$ is the surface measure of $S(\mathbb{R}^n)$. For $i > N$, we set $\mathcal{L}_i \equiv 0$. We shall always assume that Lipschitz–Killing curvatures are finite.

If you decided earlier to take the route of avoiding the generality of stratified manifolds, and you want to think of $M$ as being a smooth manifold without boundary, *and* you want to work with the standard Euclidean metric, then note that there is only one element in the outer sum of (2.2), which is the case $j = N$. This makes things much easier. Alternatively, if you want to think of $M$ as the simple, Euclidean, $N$-dimensional rectangle $\prod_{i=1}^N [0, T_i]$, then both $\widetilde{R}$ and $\widetilde{S}$ are identically zero, and it is a simple exercise to see that $\mathcal{L}_j(M) = \sum T_{i_1} \cdots T_{i_j}$, the sum being taken over the $\binom{N}{j}$ distinct choices of subscripts.

For further details, such as a proof of the fact that $\mathcal{L}_j$ are independent of the stratification; see, for example, [3, 5, 11, 23].

The fact that $\mathcal{L}_0(M)$ is equivalent to the *Euler–Poincaré characteristic* of $M$, and so independent of any Riemannian structure, is the celebrated Chern–Gauss–Bonnet theorem. The remaining Lipschitz–Killing curvatures also appear under a variety of other names, such as quermassintegrales, Minkowski, Dehn and Steiner functionals, and intrinsic volumes, although in many of these cases the ordering and normalizations are different from ours.

We now almost have enough to decipher the meaning of the Lipschitz–Killing curvatures appearing in the main result (1.3). What remains is to define the Riemannian metric under which they will be computed for most of this paper.

Recall that Riemannian metrics are defined on pairs of vectors in the tangent spaces $T_t$. Thus, taking $X_t, Y_t \in T_t M$, letting $y_i : M \to \mathbb{R}$ to be one of the components of our random function $y$, and writing $X_t y$ for the derivative of $y$ at $t$ and in direction $X_t$, we can define the (*Riemannian*) *metric induced by* $y$ as

$$(2.3) \quad g(X, Y) \equiv g_t(X_t, Y_t) \triangleq \mathbb{E}\{(X_t y_i) \cdot (Y_t y_i)\} = X_s Y_t C(s, t)|_{s=t},$$

where $C(s, t) = \mathbb{E}\{y_i(s) y_i(t)\}$ is the common covariance function of the $y_i$. (We implicitly assume that $C$ is such that $g$ is truly a Riemannian metric, in that $g_t$ is nondegenerate for each $t \in M$.)

Note from (2.3) that all the tools of Riemannian manifolds, connections, curvatures, etc., can be expressed in terms of covariances. In particular, it turns out that all of these tools also have interpretations in terms of conditional means and variances. See [3] for details.

One case, however, is worth noting now. If $y$ is locally isotropic, with unit variance derivatives in all directions, in the sense that

$$\mathbb{E}\{X_t y_i(t) Y_t y_i(t)\} = X_s Y_t C(s,t)|_{s=t} = \langle X_t, Y_t \rangle,$$

where the inner product is the usual one, then the Riemannian metric is the usual Euclidean one.

**3. Gaussian Minkowski functionals.** Having determined the meaning of the $\mathcal{L}_j(M)$ of the main result (1.3), we now turn to the $\mathcal{M}_j(D)$. Again, a tube formula approach is the easiest to take, at least initially. For this, let $\gamma \equiv \gamma_k$ be Gauss measure on $\mathbb{R}^k$, so that for Borel $A \subset \mathbb{R}^k$,

$$\gamma_k(A) = (2\pi)^{-k/2} \int_A e^{-|x|^2/2} \, dx.$$

Once again, take $D \subset \mathbb{R}^k$ to be a stratified manifold, satisfying the same conditions we set out for $M$. We also need an an additional integrability condition, a weak form of which is given in Chapter 15 of [1], but which we shall bypass by assuming that all (Euclidean) curvature tensors and second fundamental forms of $D$ and its strata are uniformly bounded.

Then there exist numbers $\mathcal{M}_j^\gamma(M) = \mathcal{M}_j^{\gamma_k}(M)$, $j \geq 1$, called the Gaussian Minkowski functionals so that, for $\rho$ small enough,

$$(3.1) \qquad \gamma_k(\text{Tube}(D,\rho)) = \gamma_k(D) + \sum_{j=1}^{\infty} \frac{\rho^j}{j!} \mathcal{M}_j^{\gamma_k}(D).$$

This gives an implicit definition of the $\mathcal{M}_j$, akin to that which Weyl's tube formula (2.1) gave for the $\mathcal{L}_j$. We shall give a more constructive definition, in the spirit of (2.2), in Section 6.1.

In many cases, the $\mathcal{M}_j^\gamma$ are quite easy to compute. For example, if $M = [u, \infty) \subset \mathbb{R}$, then since $\text{Tube}([u, \infty), \rho) = [u - \rho, \infty)$, comparing a Taylor series expansion of $\gamma(\text{Tube}(M, \rho))$ with (3.1) easily gives, for $j \geq 1$,

$$\mathcal{M}_j^{\gamma_1}([u, \infty)) = (2\pi)^{-1/2} H_{j-1}(u) e^{-u^2/2},$$

where $H_n$ is the $n$th Hermite polynomial. Additional examples, also accessible from simple calculus, can be found in [3, 20]. For a proof of (3.1), see either [20] where, to the best of our knowledge, it first appeared for manifolds with $C^2$ boundary, or Chapter 10 of [3], which treats stratified manifolds.

**4. The main result.** We now have all that is needed for understanding the result (1.3). Indeed, with just one more condition, we can state it properly.

Writing for the moment $y$ to represent one of the i.i.d. real valued components of our $\mathbb{R}^k$ valued random process $y$, assume that it is Gaussian, centered, of constant variance one and possessing almost surely $C^2$ sample paths. Assume furthermore that for each $t \in M$, the joint distributions of $(y, \partial y/\partial t_i, \partial^2 y/\partial t_j \, \partial t_k)_{i,j,k=1,\ldots,N}$ at $t$ are nondegenerate, and for some finite $K$, and all $s,t \in M$,

$$(4.1) \qquad \max_{i,j} |C_{ij}(t,t) + C_{ij}(s,s) - 2C_{ij}(s,t)| \leq K |\ln|t-s||^{-(1+\alpha)},$$

where $C_{ij} = \partial^4 C/\partial^2 t_i \, \partial^2 t_j$ is the covariance function of $\partial^2 y/\partial t_i \, \partial t_j$.

THEOREM 4.1. *Let $M \subset \mathbb{R}^N$ and $D \subset \mathbb{R}^k$ be stratified manifolds satisfying the conditions described in Sections 2 and 3. Let $y = (y_1, \ldots, y_k) : M \to \mathbb{R}^k$ be a vector valued random process, the components of which are independent, identically distributed, real valued, Gaussian processes satisfying the above conditions. Then*

$$(4.2) \quad \mathbb{E}\{\mathcal{L}_i(M \cap y^{-1}(D))\} = \sum_{j=0}^{\dim M - i} \begin{bmatrix} i+j \\ j \end{bmatrix} (2\pi)^{-j/2} \mathcal{L}_{i+j}(M) \mathcal{M}_j^\gamma(D),$$

*where the $\mathcal{L}_j$, $j = 0, \ldots, N$ are the Lipschitz–Killing measures on $M$ with respect to the metric induced by the $y_i$, the $\mathcal{M}_j^\gamma$ are the Gaussian Minkowski functionals on $\mathbb{R}^k$ and the combinatorial terms are given by*

$$(4.3) \quad \begin{bmatrix} n \\ k \end{bmatrix} = \frac{[n]!\omega_n}{[k]![n-k]!\omega_n \omega_{n-k}}, \qquad [n]! = n!\omega_n, \qquad \omega_n = \frac{\pi^{n/2}}{\Gamma(\frac{n}{2}+1)},$$

$\omega_n$ *being the volume of the unit ball in $\mathbb{R}^n$.*

Now that we know what all the terms above mean, it is worthwhile to look again at the structure of (1.3)/(4.2). Note that as far as the right-hand side of the equation goes, there is an elegant factoring of parameters. The covariance structure of $y$ appears only in the $\mathcal{L}_j(M)$, as does the geometry of the parameter space, while the $\mathcal{M}_j^\gamma(D)$ are computed without any reference to either of these. This is particularly useful for applications such as (1.5), since it implies that when working with $f$ of the form $f = F(y)$ the computations related to $M$ and the distributional properties of $y$ separate out completely from the properties of $F$.

We can now turn to the proof of Theorem 4.1, which will take up the rest of the paper. In fact, we shall give a full proof for only one, very specific process, on a very specific parameter space. This special process is defined on unit sphere $S(\mathbb{R}^l)$ and is known as the *canonical isotropic process*. Its components are independent, centered Gaussian processes with common covariance function

$$\mathbb{E}\{y_i(s)y_i(t)\} = \langle s, t \rangle, \tag{4.4}$$

where $\langle \cdot, \cdot \rangle$ is the usual Euclidean inner product. You can jump forward to Section 8 to see why this process is so special, and why all other cases effectively follow from it. Now, however, we shall begin proving Theorem 4.1 for $N$-dimensional $M \subset S(\mathbb{R}^l)$, $l \geq N$, and for the canonical process.

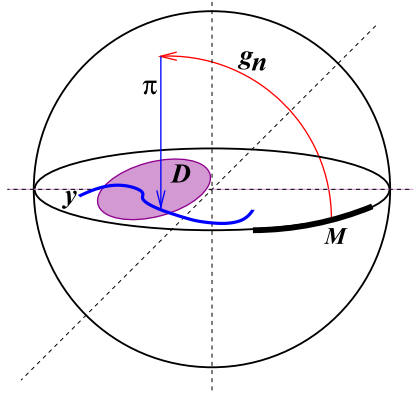

Fig. 2. *The pre-Gaussian process $y^{(2)}$ from a subset of a great circle to $\mathbb{R}^2$.*

In fact, even the canonical isotropic process will be a little too hard for us to work with directly, so we shall approximate it by something even simpler. The approximations are given by a sequence, $\{y^{(n)}\}_{n \geq l}$, of smooth $\mathbb{R}^k$-valued processes on $S(\mathbb{R}^l)$. To define them, for each $n \geq l$, we first embed $S(\mathbb{R}^l)$ in $S(\mathbb{R}^n)$ in the natural way, by setting

$$S(\mathbb{R}^l) = \{t = (t_1, \ldots, t_n) \in S(\mathbb{R}^n) : t_{l+1} = \cdots = t_n = 0\}.$$

Taking the rotation group $O(n)$, equipped with its normalized Haar measure $\mu_n$, as an underlying probability space, the $n$th process $y^{(n)}$ is defined by

$$y^{(n)}(t, g_n) \stackrel{\Delta}{=} \pi_{\sqrt{n},n,k}(\sqrt{n}g_n t), \tag{4.5}$$

where $t \in S(\mathbb{R}^l)$, $g_n \in O(n)$ and $\pi_{\sqrt{n},n,k}$ is the projection from $S_{\sqrt{n}}(\mathbb{R}^n)$ to $\mathbb{R}^k$ given by

$$\pi_{\lambda,n,k}(x_1, \ldots, x_n) = (x_1, \ldots, x_k). \tag{4.6}$$

Figure 2 gives an example of this when $M$ is a subset of a $N = 1$ dimensional great circle in the 2-sphere ($l = 2$), taking values in the plane ($k = 2$).

Now recall the old and quite elementary result known as the Poincaré limit theorem, or Poincaré's limit (although whether it really is due to Poincaré is not clear [7]). In our current setting, a recent more powerful version due to Diaconis, Eaton and Lauritzen [8] implies that the finite dimensional distributions of $y^{(n)}$ converge, in total variation norm, to those of the canonical process on $\mathbb{R}^l$. More importantly, for functionals $F$ of these processes for which $\mathbb{E}\{|F(y)|\} < \infty$, we have

(4.7) $$\lim_{n\to\infty} \mathbb{E}\{F(y^{(n)})\} = \mathbb{E}\{F(y)\}.$$

Thus, if we can compute mean Lipschitz–Killing curvatures for the approximating processes, we shall be in good shape obtaining the result we are seeking by a limiting argument. The key to the first step lies in a version of the KFF on spheres.

**5. The KFF on $S_\lambda(\mathbb{R}^n)$ and the beginnings of a proof.** Although we met the KFF for Euclidean spaces back in the Introduction [cf. (1.6)], we need a version for the sphere $S(\mathbb{R}^n)$, good references for which include [9, 13, 17, 18]. To formulate it, we need to extend the Lipschitz–Killing curvatures of Section 2 to a one parameter family of measures $\mathcal{L}_i^\kappa$, for $\kappa \geq 0$ and $0 \leq i \leq N = \dim(M)$, by replacing the Riemannian metric $\widetilde{R}$ in (2.2) by $\widetilde{R} + (\kappa/2)I^2$, where $I$ is the identity tensor. As for the basic Lipschitz–Killing curvatures, for $i > N$, we define $\mathcal{L}_i^\kappa \equiv 0$.

We also define a one parameter family of curvature integrals $\mathcal{L}_i^\kappa(M) \triangleq \mathcal{L}_i^\kappa(M, M)$. For lack of better terminology, we call these new curvature measures *extended* Lipschitz–Killing curvatures. It is trivial that $\mathcal{L}_i^0 \equiv \mathcal{L}_i$. It takes some algebra (cf. [3]), but it is not too hard to show that there are simple equivalences between the $\mathcal{L}_j$ and $\mathcal{L}_j^\kappa$, given by

(5.1) $$\mathcal{L}_i^\kappa = \sum_{n=0}^\infty \frac{(-\kappa)^n}{(4\pi)^n} \frac{(i+2n)!}{n!i!} \mathcal{L}_{i+2n}, \qquad \mathcal{L}_i = \sum_{n=0}^\infty \frac{\kappa^n}{(4\pi)^n} \frac{(i+2n)!}{n!i!} \mathcal{L}_{i+2n}^\kappa.$$

[Note that the summations here are actually finite, since all curvatures are zero for $i > \dim(M)$.]

It turns out that in dealing with subsets of $S_\lambda(\mathbb{R}^n)$ the $\mathcal{L}_j^{\lambda^{-2}}$ are more natural to deal with than are the original Lipschitz–Killing curvatures. For example, there is a elegant version of Weyl's tube formula involving them, and, more importantly for us, there is a nice KFF.

Let $G_{n,\lambda}$ denote the group of isometries (i.e., rotations) on $S_\lambda(\mathbb{R}^n)$, with Haar measure $\nu_{n,\lambda}$ normalized so that for any $x \in S_\lambda(\mathbb{R}^n)$ and every Borel $A \subset S_\lambda(\mathbb{R}^n)$, $\nu_{n,\lambda}(\{g_n \in G_{n,\lambda}: g_n x \in A\}) = \mathcal{H}_{n-1}(A)$, where $\mathcal{H}_{n-1}$ is surface (Hausdorff) measure. The KFF on $S_\lambda(\mathbb{R}^n)$ then reads as follows, where $M_1$ and $M_2$ are stratified manifolds in $S_\lambda(\mathbb{R}^n)$ satisfying the conditions that we have been assuming:

$$\int_{G_{n,\lambda}} \mathcal{L}_i^{\lambda^{-2}}(M_1 \cap g_n M_2) \, d\nu_{n,\lambda}(g_n)$$

$$
\begin{align}
(5.2) \quad &= \sum_{j=0}^{n-1-i} \begin{bmatrix} i+j \\ i \end{bmatrix} \begin{bmatrix} n-1 \\ j \end{bmatrix}^{-1} \mathcal{L}_{i+j}^{\lambda-2}(M_1) \mathcal{L}_{n-1-j}^{\lambda-2}(M_2) \\
&= \sum_{j=0}^{n-1-i} \frac{s_{i+1} s_n}{s_{i+j+1} s_{n-j}} \mathcal{L}_{i+j}^{\lambda-2}(M_1) \mathcal{L}_{n-1-j}^{\lambda-2}(M_2).
\end{align}
$$

Now let us apply this KFF to the processes $y^{(n)}$. The key result, from which everything else follows, is the following, under the usual conditions on $M$ and $D$.

LEMMA 5.1. *Let $y^{(n)}$ be the model process (4.5) on $M \subset S(\mathbb{R}^l)$, with $n \geq l$. Then for $D \subset \mathbb{R}^k$,*

$$\mathbb{E}\{\mathcal{L}_i^1(M \cap (y^{(n)})^{-1} D)\}$$

$$
\begin{align}
(5.3) \quad &= \sum_{j=0}^{\dim M - i} \left( n^{j/2} \begin{bmatrix} n-1 \\ j \end{bmatrix}^{-1} \right) \begin{bmatrix} i+j \\ j \end{bmatrix} \mathcal{L}_{j+i}^1(M) \frac{\mathcal{L}_{n-1-j}^{n-1}(\pi_{\sqrt{n},n,k}^{-1} D)}{s_n n^{(n-1)/2}} \\
&= \sum_{j=0}^{\dim M - i} \frac{s_{i+1}}{s_{i+j+1}} \mathcal{L}_{j+i}^1(M) \frac{\mathcal{L}_{n-1-j}^{n-1}(\pi_{\sqrt{n},n,k}^{-1} D)}{s_{n-j} n^{(n-1-j)/2}}.
\end{align}
$$

[It is important to understand the meaning of $\pi_{\sqrt{n},n,k}^{-1} D$ above, and in all that follows. The problem is that for all $t \in S_{\sqrt{n}}(\mathbb{R}^n)$, $\pi_{\sqrt{n},n,k}(t) \in B_{\sqrt{n}}(\mathbb{R}^n)$, which may or may not cover $D$. Thus, since

$$\pi_{\sqrt{n},n,k}^{-1} D = \{ t \in S_{\sqrt{n}}(\mathbb{R}^n) : \pi_{\sqrt{n},n,k}(t) \in D \},$$

it follows that $\pi_{\sqrt{n},n,k}^{-1} D$ may be only the inverse image of a subset of $D$.]

PROOF OF LEMMA 5.1. Since $\pi_{\sqrt{n},n,k}^{-1} D$ is a $C^2$ domain in $S(\mathbb{R}^n)$, it follows from the construction of $y^{(n)}$ that

$$\mathbb{E}\{\mathcal{L}_i^1(M \cap (y^{(n)})^{-1} D)\}$$

$$= \int_{O(n)} \mathcal{L}_i^1(M \cap (y^{(n)})^{-1} D)(g_n) \, d\mu_n(g_n)$$

$$= \int_{O(n)} \mathcal{L}_i^1(M \cap n^{-1/2} g_n^{-1}(\pi_{\sqrt{n},n,k}^{-1} D)) \, d\mu_n(g_n)$$

$$= n^{-i/2} \int_{O(n)} \mathcal{L}_i^{n-1}(\sqrt{n} M \cap g_n^{-1}(\pi_{\sqrt{n},n,k}^{-1} D)) \, d\mu_n(g_n)$$

$$= \frac{1}{s_n n^{(n-1+i)/2}} \int_{G_{n,\sqrt{n}}} \mathcal{L}_i^{n-1}(\sqrt{n} M \cap g_n(\pi_{\sqrt{n},n,k}^{-1} D)) \, d\nu_{n,\sqrt{n}}(g_n),$$

where the second-to-last line follows from the scaling properties of Lipschitz–Killing curvatures and the last is really no more than a notational change, using the definition of $\nu_{n,\lambda}$.

However, applying the KFF (5.2) to the last line above, we immediately have that it is equal to

$$\sum_{j=0}^{\dim M - i} n^{j/2} \begin{bmatrix} i+j \\ i \end{bmatrix} \begin{bmatrix} n-1 \\ j \end{bmatrix}^{-1} \frac{\mathcal{L}_{j+i}^{n-1}(\sqrt{n}M)}{n^{(i+j)/2}} \frac{\mathcal{L}_{n-1-j}^{n-1}(\pi_{\sqrt{n},n,k}^{-1} D)}{s_n n^{(n-1)/2}}$$

$$= \sum_{j=0}^{\dim M - i} n^{j/2} \begin{bmatrix} n-1 \\ j \end{bmatrix}^{-1} \begin{bmatrix} i+j \\ j \end{bmatrix} \mathcal{L}_{j+i}^{1}(M) \frac{\mathcal{L}_{n-1-j}^{n-1}(\pi_{\sqrt{n},n,k}^{-1} D)}{s_n n^{(n-1)/2}},$$

which proves the lemma. □

Suppose we send $n \to \infty$ in (5.3), which by Poincaré's limit is effectively equivalent to replacing the model process $y^{(n)}$ with a $\mathbb{R}^k$ valued canonical Gaussian $y$. Then in order for $\mathbb{E}\{\mathcal{L}_j(M \cap y^{-1}D)\}$ to be finite for the limit process $y$, we would like to have the following limits existing for each $j < \infty$:

$$(5.4) \qquad \widetilde{\rho}_j(D) \triangleq \lim_{n \to \infty} n^{j/2} \begin{bmatrix} n-1 \\ j \end{bmatrix}^{-1} \frac{\mathcal{L}_{n-1-j}^{n-1}(\pi_{\sqrt{n},n,k}^{-1} D)}{s_n n^{(n-1)/2}}.$$

A Stirling's formula computation shows that if the limit here exists, then

$$(5.5) \qquad \widetilde{\rho}_j(D) = (2\pi)^{-j/2} [j]! \lim_{n \to \infty} \frac{\mathcal{L}_{n-1-j}^{n-1}(\pi_{\sqrt{n},n,k}^{-1} D)}{s_n n^{(n-1)/2}}.$$

Sending $n \to \infty$ in Lemma 5.1 and applying Poincaré's limit (4.7), we see that if $\mathbb{E}\{|\mathcal{L}_i^1(M \cap y^{-1}D)|\} < \infty$, then

$$(5.6) \quad \mathbb{E}\{\mathcal{L}_i^1(M \cap y^{-1}D)\} = \lim_{n \to \infty} \mathbb{E}\{\mathcal{L}_i^1(M \cap (y^{(n)})^{-1}D)\}$$
$$= \sum_{j=0}^{\dim M - i} \begin{bmatrix} i+j \\ i \end{bmatrix} \mathcal{L}_{j+i}^1(M) \widetilde{\rho}_j(D).$$

This is starting to take the form of (1.3) and (4.2). The combinatorial flag coefficients are in place, but both sides of the equation are based on the $\mathcal{L}_{j+i}^1$ curvatures rather than the $\mathcal{L}_{j+i}$, and we have yet to identify the functions $\widetilde{\rho}_j$. Note the important fact, however, that on the right-hand side of the equation we have already managed a split into product form, each factor of which depends on the underlying manifold $M$ or the set $D$, but not both.

Moving to the desired curvatures, under the assumption that the limits (5.4) exist is the next step, and the result is summarized in the following lemma, again with the usual assumptions on $M$ and $D$.

LEMMA 5.2. *Let $M \subset S(\mathbb{R}^l)$ and assume that for $0 \leq j \leq \dim(M)$, the $\widetilde{\rho}_j(D)$ of (5.5) are well defined, finite and*

$$\mathbb{E}\{|\mathcal{L}_i(M \cap y^{-1}D)|\} < \infty \tag{5.7}$$

*for $y$ the $\mathbb{R}^k$ valued canonical isotropic Gaussian process on $S(\mathbb{R}^l)$. Then*

$$\mathbb{E}\{\mathcal{L}_i(M \cap y^{-1}D)\} = \sum_{l=0}^{\dim M - i} \begin{bmatrix} i+l \\ l \end{bmatrix} \mathcal{L}_{i+l}(M) \rho_l(D), \tag{5.8}$$

*where*

$$\rho_j(D) = \begin{cases} (j-1)! \sum_{l=0}^{\lfloor (j-1)/2 \rfloor} \dfrac{(-1)^l}{(4\pi)^l l! (j-1-2l)!} \widetilde{\rho}_{j-2l}(D), & j \geq 1, \\ \gamma_k(D), & j = 0. \end{cases}$$

PROOF. As usual, set $N = \dim(M)$. Combining (5.1) and (5.6), we have

$\mathbb{E}\{\mathcal{L}_i(M \cap y^{-1}D)\}$

$$= \mathbb{E}\left\{ \sum_{m=0}^{\lfloor (N-i)/2 \rfloor} \frac{1}{(4\pi)^m m!} \frac{(i+2m)!}{i!} \mathcal{L}^1_{i+2m}(M \cap y^{-1}D) \right\}$$

$$= \sum_{m=0}^{\lfloor (N-i)/2 \rfloor} \frac{1}{(4\pi)^m m!} \frac{(i+2m)!}{i!} \sum_{j=0}^{N-i-2m} \begin{bmatrix} i+2m+j \\ j \end{bmatrix} \mathcal{L}^1_{i+2m+j}(M) \widetilde{\rho}_j(D)$$

$$= \sum_{m=0}^{\lfloor (N-i)/2 \rfloor} \frac{1}{(4\pi)^m m!} \frac{(i+2m)!}{i!}$$

$$\times \sum_{j=0}^{N-i-2m} \begin{bmatrix} i+2m+j \\ j \end{bmatrix} \widetilde{\rho}_j(D)$$

$$\times \sum_{l=0}^{\lfloor (N-i-2m-j)/2 \rfloor} \frac{(-1)^l}{(4\pi)^l l!} \frac{(i+2m+j+2l)!}{(i+2m+j)!}$$

$$\times \mathcal{L}_{i+2m+j+2l}(M).$$

The passage from the above to (5.8) is now one of algebra, and the details of the three page calculation can be found on pages 408–410 of [3]. Since they add no insights, and since our aim in this paper is mainly to get the basic ideas across, we leave them to you to check. □

What remains now to do is to check that the limits (5.5) are indeed well defined, and to evaluate them. Before we can do this, we need to make a small excursion into differential geometry.

**6. A geometric interlude.** We need to better understand Gaussian Minkowski functionals and, in particular, their relationship to Lipschitz–Killing curvatures. In the following section, we shall use this information to evaluate the limits (5.5), and via them the all important functions $\rho_j(D)$ of (5.9).

6.1. *Gaussian Minkowski functionals.* We now give an explicit construction of Gaussian Minkowski functionals, rather than relying on their implicit definition via the tube formula of (3.1). First, however, we need to somewhat extend the notion of Lipschitz–Killing curvatures.

We shall continue to assume that $M$ is a $N$-dimensional stratified manifold in $\mathbb{R}^l$, locally convex and with uniformly bounded curvature. Take Borel $A \subset \mathbb{R}^l$ and $B \subset S(\mathbb{R}^l)$, retain the notation of (2.2) and define for $0 \leq i \leq l-1$, a family of *generalized Lipschitz–Killing curvature measures* supported on $M \times S(\mathbb{R}^l)$ by

$$\widetilde{\mathcal{L}}_i(M, A \times B)$$

(6.1)
$$\triangleq \sum_{j=i}^{\dim M} (2\pi)^{-(j-i)/2}$$
$$\times \sum_{m=0}^{\lfloor (j-i)/2 \rfloor} \frac{(-1)^m C(l-j, j-i-2m)}{m!(j-i-2m)!}$$
$$\times \int_{\partial_j M \cap A} \int_{S(T^\perp \partial_j M) \cap B} \mathrm{Tr}^{T_t \partial_j M}(\widetilde{R}^m \widetilde{S}_{\nu_{l-j}}^{j-i-2m})$$
$$\times \mathbb{1}_{N_t M} \mathcal{H}_{l-j-1}(d\nu_{l-j}) \mathcal{H}_j(dt).$$

For $i = l$, we define $\widetilde{\mathcal{L}}_i$ only on sets of the form $A \times S(\mathbb{R}^l)$ by setting $\widetilde{\mathcal{L}}_l(A \times S(\mathbb{R}^l)) = \mathcal{H}_l(A)$. For Borel $f : \mathbb{R}^l \times S(\mathbb{R}^l) \to \mathbb{R}$ let $\widetilde{\mathcal{L}}_i(M, f)$ denote the integral of $f$ with respect to $\widetilde{\mathcal{L}}_i$.

A change of numbering and normalization now defines the *generalized Minkowski curvature measures* as

(6.2) $$\widetilde{\mathcal{M}}_j(M, A \times B) \triangleq (j!\omega_j)\widetilde{\mathcal{L}}_{l-j}(M, A \times B).$$

With these definitions, we can now give a direct definition of the Gaussian Minkowski functionals appearing in the tube formula (3.1) by setting

(6.3) $$\mathcal{M}_j^\gamma(M) \triangleq (2\pi)^{-l/2} \sum_{m=0}^{j-1} \binom{j-1}{m} \widetilde{\mathcal{M}}_{m+1}(M, H_{j-1-m}(\langle t, \eta \rangle) e^{-|t|^2/2}).$$

Here, $H_n$ is the $n$th Hermite polynomial and $\langle t, \eta \rangle$ is the standard inner product between $t \in M$ and $\eta$ a vector in $S(\mathbb{R}^l)$.

6.2. *Warped products of Riemannian manifolds.* Note first that the set, $\pi^{-1}_{\sqrt{n},n,k}D \in S_{\sqrt{n}}(\mathbb{R}^n)$, so crucial to the arguments of Section 5 and crucial for all that will follow, can be stratified topologically into the disjoint union of the $(n-k-1+j)$-dimensional strata

(6.4) $$\widetilde{D}_{n-k-1+j} \cong D_j \times S(\mathbb{R}^{n-k}).$$

Furthermore, each such $\widetilde{D}_{n-k-1+j}$ can be written, again topologically as a disjoint union

(6.5) $$\widetilde{D}_{n-k-1+j} \simeq (D_j \cap S_{\sqrt{n}}(\mathbb{R}^k)) \sqcup (D_j \cap (B_{\mathbb{R}^k}(0,\sqrt{n}))^\circ \times S(\mathbb{R}^{n-k})).$$

(See the example in Figure 3.)

Since each $D_j \cap S_{\sqrt{n}}(\mathbb{R}^k)$ is a stratified subset of $S_{\sqrt{n}}(\mathbb{R}^k)$, its Lipschitz–Killing curvatures can be computed using the tools we have developed so far. The second set in the union is, however, somewhat more complex. Although we have written it as a product set, this is only correct up to topological equivalence. We, however, need precise Lipschitz–Killing curvatures, which with the exception of $\mathcal{L}_0$, are not topological invariants. To handle this, we need we need to break the Riemannian structure of products into a product of structures. Second, we need to treat each such product, in which each copy of $S(\mathbb{R}^{n-k})$ is likely to have a different radius, as a subset of a *Riemannian warped product*.

Recall that the Riemannian warped product of two Riemannian manifolds $(M_1, g_1)$ and $(M_2, g_2)$ with a smooth warp function $\sigma^2 : M_1 \to [0, +\infty)$ is the

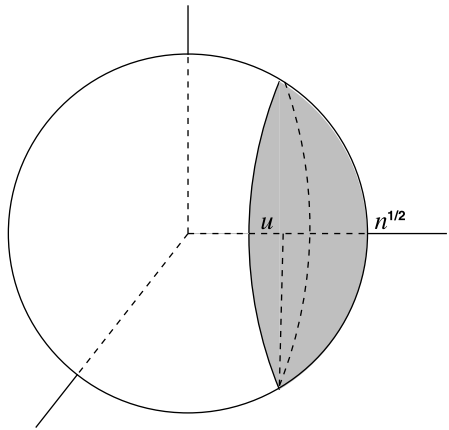

FIG. 3. *A particularly simple example of (6.5) when $N = k = 1$, $n = 3$, and $D = [u, \infty) = D_0 \sqcup D_1 \triangleq \{u\} \sqcup (u, \infty)$. In this case, $\pi^{-1}_{\sqrt{3},3,1}D \in S_{\sqrt{3}}(\mathbb{R}^3)$ is the spherical cap shown. For $j = 0$, $\widetilde{D}_{3-1-1+0} = \widetilde{D}_1$ is the set made up of the left-hand boundary of the cap along with its rightmost point. For $j = 1$, $\widetilde{D}_2$ is the rest of the cap.*

Riemannian manifold

(6.6) $$(M_1, M_2, \sigma) \triangleq (M_1 \times M_2, g_1 + \sigma^2 g_2).$$

Usually, as in our case, $M_2$ is a sphere. As an example, consider

$$\widetilde{M}_\sigma = (B_{\mathbb{R}^k}(0, \sqrt{n}))^\circ \times S(\mathbb{R}^{n-k}),$$

where the Riemannian metric on the open ball is given by

(6.7) $$g_\sigma = g_{\mathbb{R}^k} + \nabla \sigma \otimes \nabla \sigma$$

for $\sigma^2(t) = n - \|t\|_{\mathbb{R}^k}^2$, and the Riemannian metric on $S(\mathbb{R}^{n-k})$ is the canonical one inherited from $\mathbb{R}^{n-k}$ warped by $\sigma^2(t)$. The importance of this example to us is that each stratum $\widetilde{D}_{n-k-1+j}$ is isometrically embedded in such a warped product. Using this embedding, we shall be able to compute

$$\mathcal{L}_{n-1-i}^{1/n}(\pi_{\sqrt{n},n,k}^{-1} D, \widetilde{D}_{n-k-1+j}),$$

the contribution of these strata to the Lipschitz–Killing curvatures of $\pi_{\sqrt{n},n,k}^{-1} D$.

6.3. *Connections and second fundamental forms.* The first step to computing these contributions is to determine the form of the Levi–Civita connection $\widetilde{\nabla}^\sigma$ of $\widetilde{M}_\sigma$, as this is needed to in order to compute the second fundamental form of $D_j \times S(\mathbb{R}^{n-k})$ in $\widetilde{M}_\sigma$. Therefore, consider a general warped product $(M_1, M_2, \sigma)$ and denote the Levi–Civita connection on each $M_j$ by $\nabla^j$. Use $E$, or $E_j$, to denote vector fields on $M_1$, identified with their natural extensions on $M_1 \times M_2$. Similarly, $F$ or $F_j$, denote vector fields on $M_2$ extended to $M_1 \times M_2$. Then from the definition of Levi–Civita connections, the product structure of $M$ and Koszul's formula is not too hard to check that

(6.8) $$\widetilde{\nabla}^\sigma_{E_1} E_2 = \nabla^1_{E_1} E_2, \qquad \widetilde{\nabla}^\sigma_{F_1} F_2 = \nabla^2_{F_1} F_2 - \sigma^2 g_2(F_1, F_2) \nabla \sigma^2,$$
$$\widetilde{\nabla}^\sigma_E F = \widetilde{\nabla}^\sigma_F E = E(\log \sigma) F.$$

With the Levi–Civita connection defined, the next step lies in determining the second fundamental form over the sets $D_j \times S(\mathbb{R}^k)$ as they sit in $\widetilde{M}_\sigma$, as well as traces of their powers. For this, we need to describe the normal spaces $T_{(t,\eta)} \widetilde{M}_\sigma^\perp$ for $(t, \eta) \in D_j \times S(\mathbb{R}^{n-k})$. A simple argument shows that at these points,

(6.9) $$(T_{(t,\eta)} \widetilde{M}_\sigma)^\perp = (T_t D_j \oplus T_\eta S(\mathbb{R}^{n-k}))^\perp \simeq T_t D_j^\perp,$$

where $T_t D_j^\perp$ is the orthogonal (with respect to $g_\sigma$) complement of $T_t D_j$ in $T_t B_{\mathbb{R}^k}(0, \sqrt{n})$.

We can therefore state:

LEMMA 6.1. *Retaining the above notation, for $0 \leq l \leq n-1$, take*

$$(t, \eta) \in D_j \times S(\mathbb{R}^{n-k}), \qquad \nu_{k-j} \in (T_{(t,\eta)} D_j \times S(\mathbb{R}^{n-k}))^\perp.$$

*Then*

$$\frac{1}{l!} \mathrm{Tr}(S^l_{\nu_{k-j}}) = \sum_{r=0}^{l} \binom{n-k-1}{l-r} (-1)^{l-r} (\nu_{k-j}(\log \sigma_t))^{l-r} \mathrm{Tr}(S^r_{\sigma, \nu_{k-j}}),$$

*where $S$ is the second fundamental form of $D_j \times S(\mathbb{R}^{n-k})$ in $\widetilde{M}_\sigma$ and $S_\sigma$ is the second fundamental form of $D_j$ in $(B_{\mathbb{R}^k}(0, \sqrt{n}), g_\sigma)$.*

PROOF. Fix an orthonormal (under $g$) basis $(E_1, \ldots, E_k, F_1, \ldots, F_{n-k-1})$ of $T_{(t,\eta)} \widetilde{M}_\sigma$ such that $(E_1, \ldots, E_j)$ forms an orthonormal basis of $T_t D_j$. Observation (6.9) implies that any $\nu_{k-j}$ can be expressed as

$$\nu_{k-j} = \sum_{r=1}^{n-k-1} a_r F_r$$

for some constants $a_r$. Applying (6.8) and the Weingarten equation, we find

$$S_{\nu_{k-j}}(E_r, E_s) = -g(\nabla^\sigma_{E_r} \nu_{k-j}, E_s) = -g_\sigma(\nabla^\sigma_{E_r} \nu_{k-j}, E_s) = S_{\sigma, \nu_{k-j}}(E_r, E_s),$$

$$S_{\nu_{k-j}}(F_r, F_s) = -g(\nabla^\sigma_{F_r} \nu_{k-j}, F_s) = -\nu_{k-j}(\log \sigma_t) g(F_r, F_s)$$
$$= -\nu_{k-j}(\log \sigma_t) \delta_{rs},$$

$$S_{\nu_{k-j}}(E_r, F_s) = 0.$$

Therefore, for each $\nu_{k-j}$, the matrix of the shape operator $S_{\nu_{k-j}}$ in our chosen orthonormal basis is block diagonal with one block, of size $j$, given by $\{S_{\sigma, \nu_{k-j}}(E_r, E_s)\}_{1 \leq r,s \leq j}$ and the other, of size $n-k-1$, being given by $-\nu_{k-j}(\log \sigma_t) I_{(n-k-1) \times (n-k-1)}$. Therefore, applying basic combinatorial properties of the trace operator we have that for $l \leq n-k-1$,

$$\frac{1}{l!} \mathrm{Tr}(S^l_{\nu_{k-j}}) = \sum_{r=0}^{l} \binom{n-k-1}{l-r} (-1)^{l-r} (\nu_{k-j}(\log \sigma_t))^{l-r} \frac{1}{r!} \mathrm{Tr}(S^r_{\sigma, \nu_{k-j}}),$$

which completes the proof. □

**7. Back to the main proof for the special case.** With the geometry behind us, we now turn to the asymptotics required for computing the limits (5.5), for which the following lemma is a crucial step. The final move to the Gaussian Minkowski functionals appearing in Theorem 4.1 above will then involve no more than some careful asymptotics. We start with a lemma, assuming throughout that the conditions of Theorem 4.1 hold.

LEMMA 7.1. *Let $\widetilde{D}_{n-k-1+j}$ be as defined at (6.4) and $i \geq k-j \geq 0$. Then*

$$\mathcal{L}_{n-1-i}^{1/n}(\pi_{\sqrt{n},n,k}^{-1}D, \widetilde{D}_{n-k-1+j})$$

$$(7.1) \quad = s_{n-k} \sum_{l=0}^{i+j-k} \frac{s_{k+l-j}}{s_i} \binom{n-k-1}{i+j-k-l}$$

$$\times \widetilde{\mathcal{L}}_{j-l}(D^\sigma, \sigma^{n+k-2i-2j+2l-1}(2\pi)^{-k/2} h_{i+j-k-l} \mathbb{1}_{D_j}),$$

*where $D^\sigma = D \cap B_{\mathbb{R}^k}(0, \sqrt{n})$ is the regular stratified manifold obtained from the intersection of the embedding of $D$ in $\mathbb{R}^k$ and the open ball of radius $\sqrt{n}$, endowed with the metric $g_\sigma$ given by (6.7) and $h_l(t, \nu) \triangleq \langle \nu, t \rangle_{\mathbb{R}^k}^l$.*

*Furthermore, with $\varphi_k$ denoting the k-dimensional Gaussian density,*

$$(7.2) \quad \lim_{n \to \infty} \frac{1}{s_n n^{(n-1)/2}} \mathcal{L}_{n-1-i}^{1/n}(\pi_{\sqrt{n},n,k}^{-1}D, \widetilde{D}_{n-k-1+j})$$

$$= \sum_{l=0}^{i+j-k} \frac{[k+l-j]!}{[i]!} \binom{i-1}{k+l-j-1} \widetilde{\mathcal{L}}_{j-l}(D, \varphi_k(t) h_{i+j-k-l} \mathbb{1}_{D_j}).$$

PROOF. Note first that we can assume that $n$ is large enough so that

$$(7.3) \quad \widetilde{D}_{n-k-1+j} \cap S_{\sqrt{n}}(\mathbb{R}^k) = \varnothing.$$

Suppose this were not the case. Then we can further stratify $D_j$ into two parts, one of which has no intersection with $S_{\sqrt{n}}(\mathbb{R}^k)$ and one of which is contained in $S_{\sqrt{n}}(\mathbb{R}^k)$. It is only the second part that concerns us. However, for this part, we also have that $\pi_{\sqrt{n},n,k}^{-1}D_j = D_j$, and so $\pi_{\sqrt{n},n,k}^{-1}D_j$ is a $j$-dimensional stratum of $\pi_{\sqrt{n},n,k}^{-1}D$. However, such strata only contribute to Lipschitz–Killing curvatures of order 0 to $j$. Since we are interested in (asymptotically) high-order Lipschitz–Killing curvatures, we can therefore forget these components and assume that (7.3) holds.

Thus, from the definition of the $\mathcal{L}_i^\kappa$, and for all $i \geq k-j \geq 0$, we now have

$$\mathcal{L}_{n-1-i}^{1/n}(\pi_{\sqrt{n},n,k}^{-1}D, \widetilde{D}_{n-k-1+j})$$

$$= \frac{C(k-j, i+j-k)}{(2\pi)^{(i+j-k)/2}(i+j-k)!}$$

$$\times \int_{D_j \times S(\mathbb{R}^{n-k})} \int_{S(T_{(t,\eta)}D_j \times S(\mathbb{R}^{n-k})^\perp)} \operatorname{Tr}(S_{\nu_{k-j}}^{i+j-k}) \mathcal{H}_{k-j}(d\nu_{k-j})$$

$$\times \mathbb{1}_{N_t M} \mathcal{H}_{n-1-j+k}(dt, d\eta),$$

where

$$\mathcal{H}_{n-1-j+k}(dt, d\eta) = \sigma_t^{n-1-k} \mathcal{H}_{n-1-k}(d\eta) \mathcal{H}_j(dt)$$

is the Hausdorff measure that $D_j \times S(\mathbb{R}^{n-k})$ inherits from $\widetilde{D}^\sigma$, the warped product of $(D^\sigma, g_\sigma)$ and $S(\mathbb{R}^k)$ with its usual metric and warp function $\sigma^2$ as in (6.7). Then by Lemma 6.1 and (6.9),

$$\mathcal{L}^{1/n}_{n-1-i}(\pi^{-1}_{\sqrt{n},n,k}D, \widetilde{D}_{n-k-1+j})$$
$$= C(k-j, i+j-k)(2\pi)^{-(i+j-k)/2}$$
$$\times \int_{D_j} \int_{S(\mathbb{R}^{n-k})} \int_{S(T_t D_j^\perp)} \sum_{l=0}^{i+j-k} \binom{n-1-k}{i+j-k-l} \sigma_t^{n-1-k}$$
$$\times (-1)^{i+j-k-l} (\nu_{k-j}(\log \sigma_t))^{i+j-k-l}$$
$$\times \frac{1}{l!} \mathrm{Tr}(S^l_{\sigma, \nu_{k-j}}) \mathcal{H}_{k-j}(d\nu_{k-j})$$
$$\times \mathcal{H}_{n-1-k}(d\eta) \mathcal{H}_j(dt).$$

Equation (7.1) now follows from the fact that

$$\frac{C(k-j, i+j-k)(2\pi)^{-(i+j-k)/2}}{C(k-j, l)(2\pi)^{-l/2}} = \frac{s_{k+l-j}}{s_i},$$

followed by integrating over $S(\mathbb{R}^{n-k})$ and noting that

$$\nu_{k-j}(\log \sigma_t) = -\frac{\langle \nu_{k-j}, t \rangle_{\mathbb{R}^k}}{\sigma_t^2}.$$

As for the second conclusion of the lemma, (7.2), note that

$$\lim_{n \to \infty} \frac{s_{n-k}}{s_n n^{(n-1)/2}} \binom{n-k-1}{i+j-k-l} \sigma_t^{n+k-2i-2j+2l-1}$$
$$= \lim_{n \to \infty} \frac{s_{n-k}}{s_n n^{(n-1)/2}} \binom{n-k-1}{i+j-k-l} n^{(n+k-2i-2j-2l-1)/2}$$
$$\times \left(1 - \frac{\|t\|^2}{n}\right)^{(n+k-2i-2j-2l-1)/2}$$
$$= \frac{(2\pi)^{-k/2}}{(i+j-k-l)!} e^{-\|t^2\|/2}.$$

Also,

$$\frac{s_{k+l-j}}{s_i(i+j-k-l)!} = \frac{[k+l-j]!}{[i]!} \binom{i-1}{k+l-j}.$$

Finally, dominated convergence yields (7.2) and we are done. □

Theorem 4.1, for the case of the canonical process on $S(\mathbb{R}^l)$, will now follow immediately from Theorem 5.2 and the following result, modulo the issue of the moment assumption (5.7) which appears among the conditions of Lemma 5.2, but not among those of Theorem 5.2. We shall dispose of this issue in Section 8.2 below.

THEOREM 7.2. *Suppose D satisfies the conditions of Lemma 7.1. Then in the notation of that lemma, and with $\mathcal{M}_i^\gamma(D)$ defined at (6.3),*

(7.4) $$\rho_i(D) = (2\pi)^{-i}\mathcal{M}_i^\gamma(D).$$

PROOF. We start by computing the $\widetilde{\rho}_i$. By Lemma 7.1,

$$\widetilde{\rho}_i(D) = (2\pi)^{-i/2}[i]! \sum_{j=k-i}^{k-1} \sum_{l=0}^{i+j-k} \frac{[k+l-j]!}{[i]!} \binom{i-1}{k+l-j-1}$$
$$\times \widetilde{\mathcal{L}}_{j-l}(D, \varphi_k h_{i+j-k-l} \mathbb{1}_{D_j})$$
$$= (2\pi)^{-i/2} \sum_{m=0}^{i-1} \binom{i-1}{m} \widetilde{\mathcal{M}}_{m+1}(D, \varphi_k h_{i-1-m}),$$

where the $\widetilde{M}$ are the generalized Minkowski curvature measures of (6.2).

With the $\widetilde{\rho}_i$ determined, we can now turn to the $\rho_j$. By (5.9), these are

$$\rho_i(D) = (i-1)! \sum_{l=0}^{\lfloor (i-1)/2 \rfloor} \frac{(-1)^l}{(4\pi)^l l!(i-1-2l)!} \widetilde{\rho}_{i-2l}(D)$$
$$= (i-1)! \sum_{l=0}^{\lfloor (i-1)/2 \rfloor} \sum_{m=0}^{i-2l-1} \frac{(-1)^l}{(4\pi)^l l!(i-1-2l)!} (2\pi)^{-(i-2l)/2}$$
$$\times \binom{i-2l-1}{m} \widetilde{\mathcal{M}}_{m+1}(D, \varphi_k h_{i-2l-1-m})$$
$$= (2\pi)^{-i}\mathcal{M}_i^\gamma(D),$$

on applying the definitions of Hermite polynomials and $\mathcal{M}_i^\gamma(D)$. □

**8. Beyond the canonical processes.** In this section, we shall describe two different routes to extend the proof of Theorem 4.1 from the special case of the canonical process on the sphere to the more general Gaussian processes described in the theorem. The first route will indicate that this special case actually quite trivially implies many others. The second route is much longer, but gives the full result.

8.1. *The canonical process and processes with finite expansions.* To save on notation, we shall for the moment consider only real valued processes $f$. The extension to the vector valued case is trivial.

It is well known that any continuous Gaussian process on a parameter space $M$ has a representation of the form

$$f(t) = \sum_{k=1} \xi_k \varphi_k(t), \tag{8.1}$$

in which the $\xi_k$ are independent $N(0,1)$ and the $\varphi_k$ an orthonormal basis on the reproducing kernel Hilbert space of $f$. Suppose that this expansion has $l < \infty$ terms, and that $f$ has constant unit variance. Define the mapping $\varphi \colon M \to \mathbb{R}^l$ by setting $\varphi(t) = (\varphi_1(t), \ldots, \varphi_l(t))$, and note that since $1 = \mathbb{E}\{f^2(t)\} = \sum_{j=1}^{l} \varphi_k^2(t)$, $\varphi$ actually maps $M$ to a set $\varphi(M) \subset S(\mathbb{R}^l)$. Assume that $\varphi$ is one to one, and define a new process $\widetilde{f}$ on $\varphi(M)$, by setting for $x \in \varphi(M)$, $\widetilde{f}(x) = f(\varphi^{-1}(x))$. Then

$$\mathbb{E}\{\widetilde{f}(x)\widetilde{f}(y)\} = \sum_{k=1}^{l} \varphi_k(\varphi^{-1}(x))\varphi_k(\varphi^{-1}(y)) = \langle x, y \rangle. \tag{8.2}$$

That is, $\widetilde{f}$ is a version of the canonical Gaussian process on $S(\mathbb{R}^l)$.

Now note that the excursion sets of $f$ and $\widetilde{f}$ are related by the fact that $A(\widetilde{f}, \varphi(M), D) = \varphi(A(f, M, D))$ and so, as long as $\varphi$ is smooth enough ($C^2$ suffices) it is easy to relate the Lipschitz–Killing curvatures of $A(f, M, D)$ to those of $A(\widetilde{f}, \varphi(M), D)$. Their Euler characteristics, $\mathcal{L}_0$, for example, will be identical. As far as the others are concerned, the Lipschitz–Killing curvatures of a set $A \in \varphi(M)$, computed with respect to the usual Euclidean metric on $S(\mathbb{R}^l)$, will be identical to those of $\varphi^{-1}(A) \in M$ computed with respect to the Riemannian metric on $M$ which is the pull-back to $M$ of the Euclidean metric on $\varphi(M)$ by $\varphi^{-1}$. Another calculation like (8.2) shows that this is precisely the metric induced on $M$ by the process $f$.

These observations, along with the fact that Theorem 4.1 holds for the canonical process on $S(\mathbb{R}^l)$, now suffice to establish it for all processes satisfying the assumptions of the theorem and for which the expansion (8.1) can be taken to be finite.

8.2. *Completing the proof for the general case.* There is a general technique in geometry that in the context of our problem, argues as follows: Suppose that we want to prove Theorem 4.1, in particular, (4.2). Suppose we could show that there exist functions $\widetilde{\rho}(i, j, D)$ dependent on $i$, $j$ and $D$, but neither on the distribution of $y$ nor on the topology of $M$, such that

$$\mathbb{E}\{\mathcal{L}_i(M \cap y^{-1}(D))\} = \sum_{j=0}^{\dim M - i} \mathcal{L}_{i+j}(M)\rho(i, j, D). \tag{8.3}$$

Then in order to identify the functions $\rho(i, j, D)$, we could choose a parameter space and random process that were simple enough to enable us to

GAUSSIAN KINEMATIC FORMULAE 23

compute $\mathbb{E}\{\mathcal{L}_i(M \cap y^{-1}(D))\}$ in full. Writing it in the form of (8.3) would then allow us to determine the function $\rho$, and so we would have the result in full generality. In fact, this is precisely what we have done by working with the canonical Gaussian process on the sphere.

An argument of this form is, for example, generally used to prove the kinematic fundamental formulae (1.6) and (5.2).

All that remains, therefore, for us to have a full proof of Theorem 4.1 is to prove (8.3). However, as is the case for the KFF, such representation theorems are difficult to prove, requiring considerable technicalities which go far beyond the level at which we have so far worked in this paper. The basic differential geometric tool that appears, and replaces the KFF that we have relied on so far, is the critical point theory of Morse. It is this, along with the smoothness assumption (4.1), that allows us to complete all the arguments.

By agreement with the editor, who wanted to save space, and mindful of the fact that we would like to sell a few more copies of [3], we send you to Chapter 15 there, (which, in turn, will send you back to Chapter 13, which will demand that you also read Chapters 8–12) to see how this argument plays out in detail.

Many other details and additional results that we left out of this paper for space reasons can also be found there.

## REFERENCES

[1] ADLER, R. J. (1981). *The Geometry of Random Fields*. Wiley, Chichester. MR611857
[2] ADLER, R. J. (2000). On excursion sets, tube formulas and maxima of random fields. *Ann. Appl. Probab.* **10** 1–74. MR1765203
[3] ADLER, R. J. and TAYLOR, J. E. (2007). *Random Fields and Geometry*. Springer, New York. MR2319516
[4] ADLER, R. J., TAYLOR, J. E. and WORSLEY, K. J. *Applications of Random Fields and Geometry: Foundations and Case Studies*. In preparation. Available at http://ie.technion.ac.il/~radler/publications.html.
[5] BERNIG, A. and BRÖCKER, L. (2002). Lipschitz–Killing invariants. *Math. Nachr.* **245** 5–25. MR1936341
[6] BRÖCKER, L. and KUPPE, M. (2000). Integral geometry of tame sets. *Geom. Dedicata* **82** 285–323. MR1789065
[7] DIACONIS, P. and FREEDMAN, D. (1987). A dozen de Finetti-style results in search of a theory. *Ann. Inst. H. Poincaré Probab. Statist.* **23** 397–423. MR898502
[8] DIACONIS, P. W., EATON, M. L. and LAURITZEN, S. L. (1992). Finite de Finetti theorems in linear models and multivariate analysis. *Scand. J. Statist.* **19** 289–315. MR1211786
[9] FEDERER, H. (1969). *Geometric Measure Theory. Die Grundlehren der Mathematischen Wissenschaften* **153**. Springer, New York. MR0257325
[10] GORESKY, M. and MACPHERSON, R. (1988). *Stratified Morse Theory. Ergebnisse der Mathematik und Ihrer Grenzgebiete (3) [Results in Mathematics and Related Areas (3)]* **14**. Springer, Berlin. MR932724
[11] GRAY, A. (1990). *Tubes*. Addison-Wesley, Redwood City, CA. MR1044996

[12] HUG, D. and SCHNEIDER, R. (2002). Kinematic and Crofton formulae of integral geometry: Recent variants and extensions. In *Homenatge al Professor Luís Santaló i Sors* (C. BARCELÓ I VIDAL, ED.) 51–80. Univ. de Girona, Girona.
[13] KLAIN, D. A. and ROTA, G.-C. (1997). *Introduction to Geometric Probability*. Cambridge Univ. Press, Cambridge. MR1608265
[14] PFLAUM, M. J. (2001). *Analytic and Geometric Study of Stratified Spaces. Lecture Notes in Mathematics* **1768**. Springer, Berlin. MR1869601
[15] RICE, S. O. (1939). The distribution of the maxima of a random curve. *Amer. J. Math.* **61** 409–416. MR1507385
[16] RICE, S. O. (1945). Mathematical analysis of random noise. *Bell System Tech. J.* **24** 46–156. MR0011918
[17] SANTALÓ, L. A. (1976). *Integral Geometry and Geometric Probability. Encyclopedia of Mathematics and Its Applications* **1**. Addison-Wesley, Reading, MA. MR0433364
[18] SCHNEIDER, R. (1993). *Convex Bodies: The Brunn–Minkowski Theory. Encyclopedia of Mathematics and Its Applications* **44**. Cambridge Univ. Press, Cambridge. MR1216521
[19] SIEGMUND, D. O. and WORSLEY, K. J. (1995). Testing for a signal with unknown location and scale in a stationary Gaussian random field. *Ann. Statist.* **23** 608–639. MR1332585
[20] TAYLOR, J. E. (2006). A Gaussian kinematic formula. *Ann. Probab.* **34** 122–158. MR2206344
[21] TAYLOR, J. E. and ADLER, R. J. (2003). Euler characteristics for Gaussian fields on manifolds. *Ann. Probab.* **31** 533–563. MR1964940
[22] TAYLOR, J., TAKEMURA, A. and ADLER, R. J. (2005). Validity of the expected Euler characteristic heuristic. *Ann. Probab.* **33** 1362–1396. MR2150192
[23] WEYL, H. (1939). On the volume of tubes. *Amer. J. Math.* **61** 461–472. MR1507388
[24] WORSLEY, K. J. (1994). Local maxima and the expected Euler characteristic of excursion sets of $\chi^2, F$ and $t$ fields. *Adv. in Appl. Probab.* **26** 13–42. MR1260300
[25] WORSLEY, K. J. (1995). Boundary corrections for the expected Euler characteristic of excursion sets of random fields, with an application to astrophysics. *Adv. in Appl. Probab.* **27** 943–959. MR1358902
[26] WORSLEY, K. J. (1995). Estimating the number of peaks in a random field using the Hadwiger characteristic of excursion sets, with applications to medical images. *Ann. Statist.* **23** 640–669. MR1332586
[27] WORSLEY, K. J. (1997). The geometry of random images. *Chance* **9** 27–40.
[28] WORSLEY, K. J. (2001). Testing for signals with unknown location and scale in a $\chi^2$ random field, with an application to fMRI. *Adv. in Appl. Probab.* **33** 773–793. MR1875779

DEPARTMENT OF STATISTICS
STANFORD UNIVERSITY
STANFORD, CALIFORNIA 94305-4065
USA
E-MAIL: jonathan.taylor@stanford.edu
URL: www-stat.stanford.edu/~jtaylo/

INDUSTRIAL ENGINEERING AND
   MANAGEMENT
ELECTRICAL ENGINEERING
TECHNION, HAIFA
ISRAEL 32000
E-MAIL: robert@ieadler.technion.ac.il
URL: ie.technion.ac.il/Adler.phtml